\documentclass[11pt,reqno]{amsproc}

\title[]{2D Voigt Boussinesq Equations}

\author{Mihaela Ignatova}
\address{Department of Mathematics, Temple University, Philadelphia, PA 19122}
\email{ignatova@temple.edu}

\usepackage[margin=1in]{geometry}
\usepackage{amsmath, amsthm, amssymb}
\usepackage{times}
\usepackage{color}
\usepackage{hyperref}

\newcommand{\pa}{\partial}
\newcommand{\la}{\label}
\newcommand{\fr}{\frac}
\newcommand{\te}{\theta}
\newcommand{\na}{\nabla}
\newcommand{\be}{\begin{equation}}
\newcommand{\ee}{\end{equation}}
\newcommand{\ba}{\begin{array}{l}}
\newcommand{\ea}{\end{array}}
\newcommand{\R}{{\mathbb R}}
\newcommand{\TT}{{\mathbb T}}

\newcommand{\beg}{\begin}

\def\La{{\Lambda}}
\date{today}
\begin{document}
\begin{abstract} We consider a critical conservative Voigt regularization of the 2D incompressible Boussinesq system on the torus. We prove the existence and uniqueness of global smooth solutions and their convergence in the smooth regime to the Boussinesq solution when the regularizations are removed. We also consider a range of mixed (subcritical-supercritical) Voigt regularizations for which we prove the existence of global smooth solutions.  
\end{abstract}
\keywords{Boussinesq, Voigt regularization, global existence}

\noindent\thanks{\em{ MSC Classification:  35Q30, 35Q35, 35Q92.}}

\maketitle

\section{Introduction}

The Boussinesq equations are basic models of incompressible fluids in which density variations are due to variations of the temperature \cite{chandra}. They form the basis of studies of thermally generated turbulence, with applications to geophysics and theoretical physics. In addition, in the absence of molecular friction due to viscosity and neglecting thermal diffusivity, in two spatial dimensions, the Boussinesq equations are similar to the 3D axisymmetric incompressible Euler equations, and as such they have been extensively studied numerically and theoretically in the context of the famous problems of finite time singularities. This connection between the 2D  Boussinesq and 3D axisymmetric Euler singularities dates back at least as far as December 1991 \cite{pusi}, but it had a significant resurgence in recent years \cite{CH, CH22, tarek, tareki, tarekij}.

Fractional dissipative Boussinsesq equations generalize the original dissipative system. The equations have the form
\begin{align}
& \pa_t u + \nu\La^{\alpha}u + u \cdot \na u + \na p = Kg\theta e_2, \la{EQ01}\\
& \nabla \cdot u = 0, \la{EQ02}\\
& \pa_t \theta + \kappa\La^{\beta}\theta+ u \cdot \nabla \theta = 0 \la{EQ03}
\end{align}
where $u$ is the two dimensional divergence-free velocity, $\theta$ is the temperature, $p$ is the pressure, $\La = (-\Delta)^{\fr{1}{2}}$, $e_2$ is the unit vector pointing in the direction opposite to gravity. When $\alpha=\beta=2$, then we have the classical system, and $\nu>0$ is the kinematic viscosity, $\kappa>0$ is the thermal diffusivity, $g$ is the constant of gravitational acceleration, and $K>0$ is a constant thermal expansion coefficient. The reference constant density is taken to be 1. 

Global existence and uniqueness of regular solutions have been proved for the fractional dissipative Boussinesq system, when
$\alpha=2, \nu>0, \kappa=0$ or when $\nu=0, \beta=2, \kappa>0$ in \cite{Chae} and in \cite{HouLi} for smooth enough initial data. See also \cite{KWZ16} for persistence of regularity with less smooth initial data. Several works  \cite{danpai,LarLunTiti13} treat anisotropic viscosity cases and \cite{HKZ13} considers regularity in bounded domains.

Smaller powers of the fractional Laplacians have been successfully considered in \cite{HKR10, HKR11} with $\alpha=1, \nu>0, \kappa=0$ or $\nu=0, \beta=1, \kappa>0$. Global regularity for the 2D Boussinesq equations with $\nu\kappa>0$ and  fractional dissipation was obtained in \cite{SW} when $\alpha+\beta=1$, $\alpha>0.798103\dots$ and $\beta>0$ and $(u_0,\theta_0)\in H^s(\R^2)^2$ for $s>2$ (see also \cite{HS17}).

In this paper, we consider a Voigt regularization of the Boussinesq system. Voigt (or Kelvin-Voigt) equations have been introduced to model complex fluids with polymeric interactions \cite{osk}. Mathematically, they have been 
widely studied as Voigt regularizations in the context of incompressible fluid dynamics and magneto-hydrodynamics \cite{caoetal2004,lariostiti,lariostiti14,linshitztiti}.  Equations with Voigt regularizations have been used to obtain statistical solutions  \cite{levantramostiti,ramostiti}, in the limit of the regularization parameter tending to zero. More recently, a Voigt regularization approach was used in \cite{conpas} to tackle the problem of magnetic reconnection (topological change of magnetic field lines) and asymptotic approach to equilibrium.

As opposed to viscous regularizations, the Voigt regularizations are not dissipative, rather they are conservative. They work by regularizing the conserved energy. They have the advantage of preserving the steady states of the system, and are hence well suited for long time behavior studies, and also, they do not introduce spurious boundary layers. On the other hand, the Voigt regularizations do not respect the transport structure of the equations.

We consider the 2D Voigt Boussinesq equation,
\begin{align}
& (I + \epsilon \La)\pa_t u + u \cdot \na u + \na p = \theta e_2, \la{EQ1}\\
& \nabla \cdot u = 0, \la{EQ2}\\
& (I + \epsilon \La)\pa_t \theta + u \cdot \nabla \theta = 0 \la{EQ3}
\end{align}
on the torus $\TT^2$, where $u = (u_1, u_2) \colon \TT^2 \times [0, T] \to \R^2$ is the fluid velocity, 
$\theta \colon \TT^2 \times [0, T] \to \R$ is the temperature,
$I$ is the identity operator, $\Lambda = (- \Delta)^{\frac{1}{2}}$, $\epsilon > 0$ is a parameter, and $e_2 = (0, 1)$.
We study the initial value problem for the system \eqref{EQ1}--\eqref{EQ3} with initial data
\be
u(\cdot, 0) = u_0 \quad \text{and} \quad \theta(\cdot, 0) = \theta_0. \la{EQ4}
\ee
The curl $\omega = \nabla^\perp \cdot u$, with $\nabla^\perp = (- \partial_2, \partial_1)$ obeys the vorticity equation
\[
(I + \epsilon \La)\pa_t\omega + u\cdot\na \omega = \pa_1\theta.
\]
Thus, the Voigt Boussinesq equations \eqref{EQ1}--\eqref{EQ3} can be closed in terms of vorticity and temperature $(\omega,\theta)$,
\begin{align}
& (I + \epsilon \Lambda)\pa_t \omega + u \cdot \nabla \omega =  \partial_1 \theta, \la{EQ5}\\
& u = \nabla^\perp \Delta^{-1} \omega, \la{EQ6}\\
& (I + \epsilon \Lambda)\pa_t \theta + u \cdot \nabla \theta = 0 \la{EQ7}
\end{align}
with initial conditions 
\be
\omega(\cdot, 0) = \omega_0 = \nabla^\perp \cdot u_0 \quad \text{and} \quad \theta(\cdot, 0) = \theta_0.\la{EQ8}
\ee

Our main result is the global regularity of the Voigt Boussinesq equations.
\beg{thm}\la{globalreg}
Let $s>1$ and let $(\omega_0, \theta_0)\in (H^s(\TT))^2$. Let $T$ be arbitrary. Then, there exists a unique solution of 2D Voigt Boussinesq equations $(\omega,\te) \in L^{\infty}(0, T; (H^s(\TT^2))^2)$ with initial data 
$(\omega_0, \theta_0)$.
\end{thm}

The proof is based on the following elements. We first prove the solutions $(\omega,\te)$ exist and are unique locally in time 
in the framework of Sobolev spaces $H^s\times H^s$, as soon as $s>1$ (Theorem~\ref{local}). Then, for $s>1$, for any local solution defined on some time interval, we exploit the fundamental conservation property of Voigt regularizations. This conservation is shown to imply a priori information on the velocity  $u\in L^{\infty}(0,T, H^{\fr{3}{2}})$, with bounds that depend only on the initial data and grow at most linearly in $T$. This is used, in conjunction with the local existence and uniqueness theorem to obtain global unique solutions in the Sobolev phase space $H^s\times H^s$ for $1<s\le\fr{3}{2}$. Next, we prove a natural Beale-Kato-Majda-type theorem (Theorem~\ref{BKM}) that gives sufficient conditions for persistence of regularity which do not involve the Lipschitz norm of the temperature. Finally, based on the Theorem~\ref{BKM}, we prove that global regularity holds for any $s>\fr{3}{2}$.

In view of the interest in the blow up problem, it is important to study various relaxations of the regularization. In this paper we prove 
that the limit of vanishing regularization is the original equation, in a smooth enough regime.
\beg{thm}\la{limit} Let $(\omega_0, \theta_0)\in H^s(\TT^2)\times H^{s+1}(\TT^2)$, $s>1$ and let $(\omega_B, \theta_B)\in L^{\infty}(0,T; H^s(\TT^2)\times H^{s+1}(\TT^2))$ be a solution of the 2D Boussinesq system on $[0,T]$. Then the solutions $(\omega_{\epsilon}, \theta_{\epsilon})$ of the 2D Voigt Boussinesq equations with the same initial data converge as $\epsilon\to 0$ 
 in $L^{\infty}(0,T; H^{-1}(\TT^2)\times L^2(\TT^2))$ to the solution $(\omega_B, \te_B)$ of the 2D Boussinesq system.
\end{thm}

\beg{rem} The $s+1$ regularity requirement for $\theta$ is needed because of the conservative nature of the Boussinesq equations, to ensure $\na\theta\in L^{\infty}$ uniformly in $\epsilon$.
\end{rem}

We consider also fractional Voigt Boussinesq equations and prove global regularity for certain cases whereby lower powers of the fractional Laplacian are used for the temperature field.
\beg{thm}\la{frac} Let $s>1$. Let $\alpha, \beta\ge 0$ with $\alpha+\beta\ge 2$, $\alpha>1$,
$\beta\ge \fr{2}{3}$. Let $(\omega_0,\theta_0) \in (H^s(\mathbb T^2))^2$. Let $T>0$ be arbitrary.
Then there exists a unique solution $(\omega, \theta)\in L^{\infty}(0,T; (H^s(\mathbb T^2))^2)$ of the fractional Voigt Boussinesq system
\be
\left\{
\ba
(I + \epsilon\La)^{\alpha}\pa_t u + u\cdot\na u + \na p = \theta e_2,\\
\na\cdot u = 0,\\
(I + \epsilon\La)^{\beta}\pa_t\theta + u\cdot\na \theta = 0.
\ea
\right.
\la{fracbous}
\ee
\end{thm}
The proof of this result is based on  commutator estimates and on bounds on $\theta $ in $H^1$.  We obtain also

\beg{thm} \la{alphabig} Let $\alpha>2$, and $\beta=0$. 
 Let $(\omega_0,\theta_0) \in (H^s(\mathbb T^2))^2$. Let $T>0$ be arbitrary.
Then there exists a unique solution $(\omega, \theta)\in L^{\infty}(0,T; (H^s(\mathbb T^2))^2)$ of the fractional Voigt Boussinesq system \eqref{fracbous}
\end{thm}
\beg{rem} Theorem~\ref{frac} and Theorem~\ref{alphabig} imply regularity for $\beta\ge\fr{2}{3}$ when $\alpha\ge\fr{4}{3}$ and
for $\beta =0$ when $\alpha>2$. It would be natural to conjecture global regularity for the whole range $\alpha + \beta\ge 2$, but at present we do not know how to prove this result.
\end{rem}
\section{Proof of Theorem \ref{globalreg}} 
\subsection{Local existence and uniqueness in $H^s, s>1$}
The 2D Voigt Boussinesq equations can be written in divergence form as
\be
\left\{
\ba
\pa_t \omega = -(I + \epsilon \La)^{-1}\na\cdot (u\omega) + (I + \epsilon\La)^{-1}\pa_1\theta\\
\pa_t \theta = - (I + \epsilon\La)^{-1}\na\cdot (u\theta),
\ea
\right.
\la{vote}
\ee
where we used that $u$ is divergence-free.
\beg{thm}\la{local} Let $(\omega_0,\theta_0)\in (H^s(\TT^2))^2$ with $s>1$. There exists a time $T$ depending only on $\epsilon>0$ and the norms of $\omega_0$ and $\theta_0$ in $H^s(\TT^2)$ and a unique solution $(\omega,\te )$ of the 2D Voigt Boussinesq equations, with initial data $(\omega_0, \theta_0)$ and with  $(\omega, \te )\in L^{\infty}(0, T; (H^s(\TT^2))^2)$.
\end{thm}

\noindent{\bf{Proof of Theorem \ref{local}}}.
The nonlinearity in the 2D Voigt Boussinesq equations is 
\be
\left (
\ba 
\omega\\
\te
\ea
\right )
\mapsto
\left (
\ba 
-(I+ \epsilon \La)^{-1}\na\cdot (u\omega) +  (I+ \epsilon \La)^{-1}\pa_1\theta\\
-(I+ \epsilon \La)^{-1}\na\cdot(u\theta)
\ea
\right)
=
\left (
\ba
N_1(\omega, \te)\\
N_2(\omega, \te)
\ea
\right ).
\la{nonl}
\ee
This follows from \eqref{vote} by using that $u$ is divergence-free.
For $s>1$, we have that $H^s(\TT^2)$ is a Banach algebra, and that $ H^s(\TT^2)\subset L^{\infty}(\TT^2)$. Moreover, the map $(I+\epsilon\La)^{-1}\na$ is bounded in $H^s(\TT^2)$. We deduce
\be
\ba
\|N_1(\omega,\te)\|_{H^s} \le C \|\omega\|_{H^s}^2 + \|\te\|_{H^s}\\
\|N_2(\omega, \te)\|_{H^s} \le C \|\omega\|_{H^s}\|\te\|_{H^s}
\ea
\la{nonlb}
\ee
by using also the Poincar\'{e} inequality $\|u\|_{H^s}\le C \|\omega\|_{H^s}$.
The local existence and uniqueness of solutions follows from a fixed point argument and the fact that the nonlinearity (bilinear continuous plus linear continuous) is Lipschitz in the function space $(H^s(\TT^2))^2$:
\[||N(w_1)-N(w_2)||_{(H^s(\TT^2))^2}\le L_{B}||w_1-w_2||_{(H^s(\TT^2))^2}, \quad w_1,w_2 \in B,\]
where $w_i=(\omega_i,\theta_i)$ and $B$ is a ball in $(H^s(\TT^2))^2$.
We omit  further details.\quad $\Box$

From now on, until we discuss convergence as $\epsilon\to 0$, for simplicity of exposition we set $\epsilon =1.$

\subsection{Basic Voigt energy bounds on $H^s$, $s>1$ solutions of the 2D Voigt Boussinesq equations}
 
Multiplying the temperature equation \eqref{EQ7} by $\te$ and integrating, we obtain
\be
\|\te\|^2_{L^2} + \|\te\|_{H^{\fr{1}{2}}}^2 = \|\te_0\|^2_{L^2} + \|\te_0\|_{H^{\fr{1}{2}}}^2 = A_0.\la{EQ10}
\ee
Multiplying the vorticity equation \eqref{EQ5} by $\omega$ and integrating, we have
\be
\fr{d}{dt} (\|\omega\|^2_{L^2} + \|\omega\|_{H^{\fr{1}{2}}}^2) = 2\int_{\TT^2} (\pa_1\te) \omega dx
 \le 2\|\omega\|_{H^{\fr{1}{2}}}\sqrt{A_0}.\la{EQ11}
\ee
It follows that
\be
\fr{d}{dt} \|\omega\|_{H^{\fr{1}{2}}} \le \sqrt{A_0},
\la{omh12e}
\ee
\be
\|\omega (\cdot, t)\|_{H^{\fr{1}{2}}} \le   \|\omega_0\|_{H^{\fr{1}{2}}}  + t\sqrt{A_0},
\la{omh12b}
\ee
and thus $(\omega, \te) \in L^{\infty}(0,T; (H^{\fr{1}{2}}(\TT^2))^2)$ with a priori  bounds in terms of initial data that grow at most linearly in $T$. 

\subsection{Proof of global existence and uniqueness for $1<s\le \fr{3}{2}$}
The aim is to show that local solutions in $H^s$ for  $1<s\le \fr{3}{2}$ satisfy $H^s$ bounds
\be
\fr{d}{dt} (\|\te\|_{H^s} + \|\omega\|_{H^s}) \le C (\|u\|_{H^{s}} +1)(\|\te\|_{H^s} + \|\omega\|_{H^s}).\la{EQ12}
\ee
Then, because of the a priori bound due to the basic energy structure \eqref{omh12b}, specifically that 
$\omega\in L^{\infty}(H^{\fr{1}{2}})$, it follows that $u\in L^{\infty}(H^{\fr{3}{2}})\subset L^{\infty}(H^s)$ with bounds that depend only on initial data, and grow at most linearly in time. The inequality \eqref{EQ12} implies that the $H^s$ norms of $\omega$ and $\te$ are finite.
Thus, by the local existence and uniqueness in $H^s$, obtained in Theorem~\ref{local}, the solution can be continued indefinitely.

In order to obtain \eqref{EQ12} we use the calculus inequality
\[
\|fg\|_{H^s} \le C\left ( \|f\|_{L^{\infty}} \|g\|_{H^s} + \|f\|_{H^s}\|g\|_{L^{\infty}}\right)
\]
valid for any $s>0$, to obtain from \eqref{vote}, using the boundedness of Riesz transforms in $H^s$ spaces,
\be
\fr{d}{dt}(\|\omega\|_{H^s} +\| \te \|_{H^s})
\le \|u\|_{L^{\infty}} (\|\omega\|_{H^s} + \|\te\|_{H^s})+ (\|\omega\|_{L^{\infty}} + \|\te\|_{L^{\infty}})\|u\|_{H^s} + \|\te\|_{H^s}.\la{EQ13}
\ee
We then apply the embedding $H^s\subset L^{\infty}$ to estimate $\|\omega\|_{L^{\infty}} + \|\te\|_{L^{\infty}}$ by
$\|\te\|_{H^s} + \|\omega\|_{H^s}$ and $\|u\|_{L^{\infty}}$ by $\|u\|_{H^s}$. \quad $\Box$

\subsection{ A Beale-Kato-Majda-type theorem}
We give here a regularity criterion which can be used to uniquely continue solutions. 

\beg{thm}\la{BKM} Let $s>1$. Let $(\omega, \te)\in  L^{\infty}(0, T; (H^s(\TT^2))^2) $ be the solution of 2D Voigt Boussinesq equations with initial data $(\omega_0, \te_0)\in (H^s(\TT^2))^2$. If
\be
\int_0^T(\|\omega (\cdot, t)\|_{L^{\infty}}+ \|\te(\cdot,t)\|_{L^{\infty}})dt <\infty \la{EQ14}
\ee
then the solution $(\omega, \te )$ can be uniquely continued for a time interval $[0, T_1]$, with $T_1>T$.
In particular, $T$ is not a blow up time in $H^s$.
\end{thm}
\noindent{\bf{Proof of Theorem \ref{BKM}}}.
We return to the inequality \eqref{EQ13}. From the basic energy estimates,
$(\omega,  \theta)$
are bounded $L^{\infty}(0,T; H^{\fr{1}{2}}(\TT^2)^2)$ in terms of the initial data and $T$. The Sobolev embedding $H^{\fr{3}{2}}\subset L^{\infty}$ implies in particular that 
\[
\int_0^T \|u\|_{L^{\infty}}dt <\infty
\] 
is bounded a priori in terms of initial data and $T$, so that the only amplification factor that remains to be controlled in  \eqref{EQ13}  is $\|\omega\|_{L^{\infty}} + \|\te\|_{L^{\infty}}$, which is then controlled by the assumption \eqref{EQ14}.   We omit further details. \quad $\Box$

\beg{rem}
We note that for the Voigt-Boussinesq system the Beale-Kato-Majda criterion does not require the finiteness of $\|\nabla\theta\|_{L^1(0,T; L^{\infty})}$, as opposed to the criterion for the Boussinesq system \cite{CKN}. 
\end{rem}
\subsection{Proof of global existence and uniqueness for $s>\fr{3}{2}$}
We use the result for $s=\fr{3}{2}$ 
which gives a priori bounds for $(\omega,\te)$ in $L^{\infty}(0,T; (H^{\fr{3}{2}}(\TT^2))^2)$. 
 Then the embedding $H^{\fr{3}{2}}\subset L^{\infty}$ implies that the BKM quantity \eqref{EQ14} is controlled, and higher regularity is propagated.\quad $\Box$

\section{Proof of Theorem \ref{limit}}
Denoting by $\te_B$ and $u_B$ the solutions of the 2D Boussinesq equations, by $\te_V$ and $u_V$ the solutions of the 2D Voigt Boussinesq equation, we set
\[
u = u_V- u_B,  \quad \te = \te_V-\te_B,
\]
and we have the equations
\be
\ba
\pa_t(I +\epsilon\La)\te + (u + u_B)\cdot\na \te + u\cdot\na {\te}_B = -\epsilon\La\pa_t\te_B,\\
\pa_t(I + \epsilon\La)u+ (u+ u_B)\cdot\na u + u\cdot\na {u}_B - \theta e_2 = -\epsilon\La \pa_t u_B,\\
\na\cdot u = 0,
\ea
\ee
with vanishing initial data. We multiply the first equation by $\te$, the second equation by $u$, integrate and add.
Because $\omega_B  \in L^{\infty}(H^{s})$ and $\te_B \in L^{\infty}(H^{s+1})$,  $s>1$,   we have 
$\na {u}_B \in L^{\infty}$ and $\na{\te}_B\in L^{\infty}$ bounded uniformly by a constant $ {C}$. Hence 
\[
\left |\int_{\TT^2}( (u\cdot\na{\te}_B)\te +(u\cdot\na {u}_B) u)dx\right |
 \le {{C}}\|u\|_{L^2}(\|\te\|_{L^2}+ \|u\|_{L^2}). 
\]
We consider the quantity 
\[
E = \|u\|_{L^2}^2 + \|\theta\|_{L^2}^2 + \epsilon\left( \|u\|_{H^{\fr{1}{2}}}^2 + \|\te\||_{H^{\fr{1}{2}}}^2\right)
\]
and note that we have
\[
\fr{d}{dt} E \le C_1 (\|u\|_{L^2}^2 + \|\theta\|_{L^2}^2) + \epsilon \left[\|\pa_t\te_B\|_{H^{\fr{1}{2}}}\|\te\|_{H^{\fr{1}{2}}} + \|\pa_t u_B\|_{H^{\fr{1}{2}}}\|u\|_{H^{\fr{1}{2}}}\right]
\]
Using the Boussinesq equation, the assumptions $\omega_B\in L^{\infty}(H^s), \te_B\in L^{\infty}(H^{s+1})$, $s>1$, ensure that
$\|\pa_t\te_B\|_{H^{\fr{1}{2}}}$ and  $\|\pa_t u_B\|_{H^{\fr{1}{2}}}$ are bounded uniformly by a constant. From Young's inequality we have
\[
\fr{d}{dt} E \le C_2E + \epsilon C_3 
\]
with initial data $E(0)=0$, and thus we have the convergence  $\lim_{\epsilon\to 0}E(t) =0 $ uniformly on finite time intervals.

\section{Proof of Theorem \ref{frac}} 
We consider the 2D fractional Voigt Boussinesq 
\begin{align}
& \pa_t \omega + (I + \La)^{-\alpha}(u \cdot \nabla \omega) =  (I + \La)^{-\alpha}\partial_1 \theta, \la{EQ25}\\
& u = \nabla^\perp \Delta^{-1} \omega, \la{EQ26}\\
& \pa_t \theta + (I + \La)^{-\beta}(u \cdot \nabla \theta )= 0 \la{EQ27}
\end{align}
written in divergence form with initial conditions 
\be
\omega(\cdot, 0) = \omega_0 = \nabla^\perp \cdot u_0 \quad \text{and} \quad \theta(\cdot, 0) = \theta_0.\la{EQ28}
\ee
We have set the parameter $\epsilon=1$ for simplicity.

Multiplying equation \eqref{EQ27} by $(I + \La)^{\beta}\theta$ and integrating,
\be
\frac12\frac{d}{dt}|(I + \La)^{\beta/2}\theta|_{L^2}^2=0, \la{EQ29}
\ee
multiplying equation \eqref{EQ25} by $(I + \La)^{\alpha}\omega$ and integrating,
\be
\frac12\frac{d}{dt}|(I + \La)^{\alpha/2}\omega|_{L^2}^2
=\int \partial_1 \theta \omega
\le C ||\theta||_{H^{\beta/2}}||\omega||_{H^{\alpha/2}} \la{EQ30}
\ee
if $\alpha/2+\beta/2\ge1$. 
By \eqref{EQ29}, we get $\theta\in L^{\infty}(H^{\beta/2})$ with $||\theta||_{H^{\beta/2}}\le ||\theta_0||_{H^{\beta/2}}$. By \eqref{EQ30}, we have $\omega\in L^{\infty}(H^{\alpha/2})$. 
Now we consider $s>1$ and look at the evolution of the $H^s$ norm of $\omega$. We use the equivalent form 
$\|\omega\|_s \sim \|(I + \La)^s\omega\|_{L^2}$.  We have from \eqref{EQ30}
\be
\fr{1}{2}\fr{d}{dt}\|\omega\|_s^2 = A +B
\la{AB}
\ee
with
\be
A =  -\int (I + \La)^{s-\fr{\alpha}{2}}(u \cdot \nabla \omega)(I + \La)^{s-\fr{\alpha}{2}}\omega dx
\la{A}
\ee
 and
 \be
 B = \int  (I + \La)^{s-\alpha}\pa_1\theta (I + \La)^s\omega dx
 \la{B}
 \ee
 Now we use incompressibility, integration by parts and the commutator estimate
 \be
 \|[(I + \La)^{\sigma}, u\cdot\na] \omega\|_{L^2} \le C \|\na u\|_{L^{\infty}}\|\omega\|_{\sigma}
 \la{commu}
 \ee
valid for $\sigma>0$ \cite{KP}, together with the embedding inequality
 \be
 \|\na u\|_{L^{\infty}}\le C\|\omega\|_{s}
 \la{embed}
 \ee
 which is true because $s>1$, to conclude that
 \be
 |A| \le C \|\omega\|_s\|\omega\|_{s-\fr{\alpha}{2}}^2.
\la{Aboundone}
\ee
We choose first $s=\alpha>1$ and deduce in view of \eqref{EQ29} and \eqref{EQ30} that
\be
|A| \le C \|\omega\|_s\|\omega\|_{\fr{\alpha}{2}}^2\le (\|\omega_0\|_{\fr{\alpha}{2}} + C_{\alpha}t\|\theta_0 \|_{\fr{\beta}{2}})^2\|\omega\|_{s} = (a +bt)^2\|\omega\|_s
\la{Abound}
\ee
holds for all $t\ge 0$.  When $s=\alpha$ we see that 
\be
|B|\le C\|\theta\|_1 \|\omega\|_s
\la{Bbound}
\ee
and we have to estimate the $H^1$ norm of $\theta$. This evolves according to
\be
\fr{1}{2}\fr{d}{dt} \|\theta\|_1^2 = -\int (I+ \La)^{1-\fr{\beta}{2}}(u\cdot\na \theta) (I+ \La)^{1-\fr{\beta}{2}}\theta dx
\la{evhonetheta}
\ee
and using the commutator estimate \eqref{commu}, we obtain
\be
\fr{1}{2}\fr{d}{dt}\|\theta\|_1^2 \le C\|\omega\|_s \|\theta\|_{1-\fr{\beta}{2}}^2
.\la{h1boundone}
\ee
Interpolating we may write
\be
\|\theta\|_{1-\fr{\beta}{2}}^2 \le \|\theta\|_1\|\theta\|_{1-\beta} \le C \|\theta\|_1\|\theta_0\|_{\fr{\beta}{2}}
\la{interpol}
\ee
where we used $\beta \ge \fr{2}{3}$ and \eqref{EQ29}. 
We have thus from \eqref{h1boundone} and\eqref{interpol},
\be
\fr{1}{2}\fr{d}{dt}\|\theta\|_1^2 \le C_{\beta}\|\omega\|_s \|\theta\|_{1}
.\la{h1bound}
\ee
The ODE inequality system in $X=\|\omega\|_s$ (with $s=\alpha$) and $Y=\|\theta\|_1$,  
\be
\left \{
\ba
\fr{1}{2}\fr{d}{dt} X^2 \le (a+bt)^2 X + CXY\\
\fr{1}{2}\fr{d}{dt} Y^2 \le C_{\beta}XY
\ea
\la{ode}
\right.
\ee
follows from \eqref{AB}, \eqref{Abound}, \eqref{Bbound}, \eqref{h1bound}. From this ODE we deduce a priori bounds for the quantities $X$ and $Y$ which are finite for all $t$. Once we know these bounds, we know that $\|\na u\|_{L^{\infty}}$ is controlled.
The evolution of $\|\theta\|_{s}$ then obeys, using the commutator estimate \eqref{commu},
\be
\fr{d}{dt}\|\theta\|_{s}^2 \le C \|\na u\|_{L^{\infty}} \|\theta\|_{s-\fr{\beta}{2}}^2,
\la{hstheta}
\ee
and it implies a priori bounds on $\|\theta\|_s$. The bound \eqref{Aboundone} is replaced by
\be
|A| \le C\|\na u\|_{L^{\infty}}\|\omega\|_{s-\fr{\alpha}{2}}^2
\la{Aboundnew}
\ee
and \eqref{Bbound}, is replaced by
\be
|B|\le C\|\theta\|_s \|\omega\|_s
\la{Bboundnew}
\ee
where we used $\alpha>1$. Then from \eqref{AB} it follows that $\|\omega\|_s$ is controlled.\quad $\Box$

\section{Proof of Theorem \ref{alphabig}}
We start by noticing that because $\beta=0$ we have that
\be
\|\theta (t)\|_{L^p}\le \|\theta_0\|_{L^p}
\la{lptheta}
\ee
holds for all $t\ge 0$ and all $1\le p \le \infty$. 
Then multiplying the $u$ equation \eqref{EQ27} by $(I +\La)^{\alpha} u$ and integrating, we obtain
\be
\fr{1}{2}\fr{d}{dt} \|u\|_{\fr{\alpha}{2}}^2   = \int \theta u_2 dx\le \|\theta_0\|_{L^2} \|u\|_{L^2}
\la{enbal}
\ee
from whence we deduce that
\be
\|u(t)\|_{\fr{\alpha}{2}} \le C _ Dt.
\la{ualpha}
\ee
Let us consider first $s=1+\epsilon$ where $0<\epsilon= \alpha-2$. Multiplying the $\omega$ equation
\eqref{EQ25} by $(I + \La)^{2s}\omega$ and integrating, we have
\be
\fr{1}{2}\fr{d}{dt}\|\omega\|_s^2 = -\int (I + \La)^{s-\alpha}(\na(u\omega)) (I + \La)^s\omega dx  +
\int(I + \La)^{s-\alpha}\pa_1\theta (I + \La)^s\omega dx.
\la{hsomega}
\ee
We note that $s-\alpha +1=0$ by our choice of $s$, so we obtain that
\be
\fr{1}{2}\fr{d}{dt}\|\omega\|_s^2\le C \|u\|_{L^4}\|\omega\|_{L^4}\|\omega\|_s + C \|\theta_0\|_{L^2}\|\omega\|_{s}
\la{hsomegab}
\ee
Thus we have that $\|\omega\|_s$ is bounded on $[0,T]$. Because $s>1$, implies that $\|\na u\|_{L^{\infty}}$ is bounded on $[0,T]$ and the regularity for arbitrary $s>1$ follows as in the proof of Theorem \ref{frac} above. \quad $\Box$

\vspace{.5cm}

\noindent{\bf{Conflict of interest}} The author declares that she has no conflicts of interest.\\

\noindent{\bf{Acknowledgments.}} We acknowledge discussions with Jingyang Shu. This work was partially supported by NSF grant DMS-1713985.

\end{document}